\documentclass[12pt,a4paper]{article}
  \usepackage[utf8]{inputenc}

\usepackage{amsfonts}
\usepackage{graphicx}
\usepackage{epstopdf}
\usepackage{amsmath}
\usepackage{amsthm}
\usepackage{latexsym,bm}
\usepackage{amssymb}
\usepackage{indentfirst}
\usepackage{mathrsfs}

\usepackage{hyperref}

\usepackage{graphicx}
 
 \usepackage{epigraph}  %
 \usepackage{subfigure}
 
\newtheorem{thm}{Theorem}[section]

\newtheorem{rem}[thm]{Remark}

\numberwithin{equation}{section}

\title{  ON A  LIST OF 
  ORDINARY DIFFERENTIAL EQUATIONS PROBLEMS 
\thanks{This essay contains a very free translation with comments, updates, annotations, additions, corrections    and abridgment  of 
	 ``Uma Lista de Problemas de E.D.O", \cite{l}.}}
\author{Jorge Sotomayor 	\thanks{
	The author is fellow of  CNPq. Grant: PQ-SR- 307690/2016-4.
}}

\begin{document}
	
	\begin{titlepage}
		\maketitle
	\end{titlepage}

\begin{abstract} 
	This  evocative essay 
	focuses on the mathematical  activities witnessed by the author along 1962-64 at IMPA.  
	The list of research problems proposed  in September 1962  by Mauricio Peixoto at the Seminar on   the
	Qualitative Theory  of Differential Equations
	is pointed out as a landmark for the genesis of the research interest in the Qualitative Theory of Differential Equations and 
	Dynamical Systems in Brazil.  \\
	
	
\end{abstract}

\vskip -0,5cm
 \begin{epigraph}{{\it\small Just as   every human  un\-der\-ta\-king pursues certain objectives, so also
mathematical research re\-qui\-res its problems.   \\  
It  is by the solution of problems that the investigator  
tests the temper of his steel; he finds new methods and 
 new outlooks,  and gains a wider and freer horizon.}\\

D. Hilbert, ICM,  
Paris. 1900.}

 \end{epigraph}

\section{The diary of a journey.}\label{sec:01}
To find an old mathematics notebook,
protected from total  deterioration  by a plastic cover, was 
  like 
  an 
 encounter with
the diary of a distant journey. The contact with its  
yellowish  pages
triggered
an avalanche of recollections of   the years  $1962 - 64,$ when I was a   doctoral student  in Brazil.
Its sequential structure  
prompted   the reconstruction of  the chronology of my initiation into  mathematical research.

The essential landmarks -- the stations -- of a mathematical 
peregrination,   starting   from a bookish approach, heading 
toward  an  attempt to tackle
research
problems, 
pulsed latent in the rough  writing and scribbled drawings. Some pages were missing, a few of them
 had
faded away.

\section {An   afternoon in
September 1962. } 
 
In a seminar room of  the Institute of Pure and Applied Mathematics  (IMPA),
 in the Botafogo quarter  of the city of Rio de Janeiro,  gathered  a group of around ten people:
mathematicians,
research fellows of disparate backgrounds, candidates to become mathematicians, and one or two  {\it voyeurs.}

\begin{figure}[!ht]
\begin{center}
\includegraphics[scale=0.37]{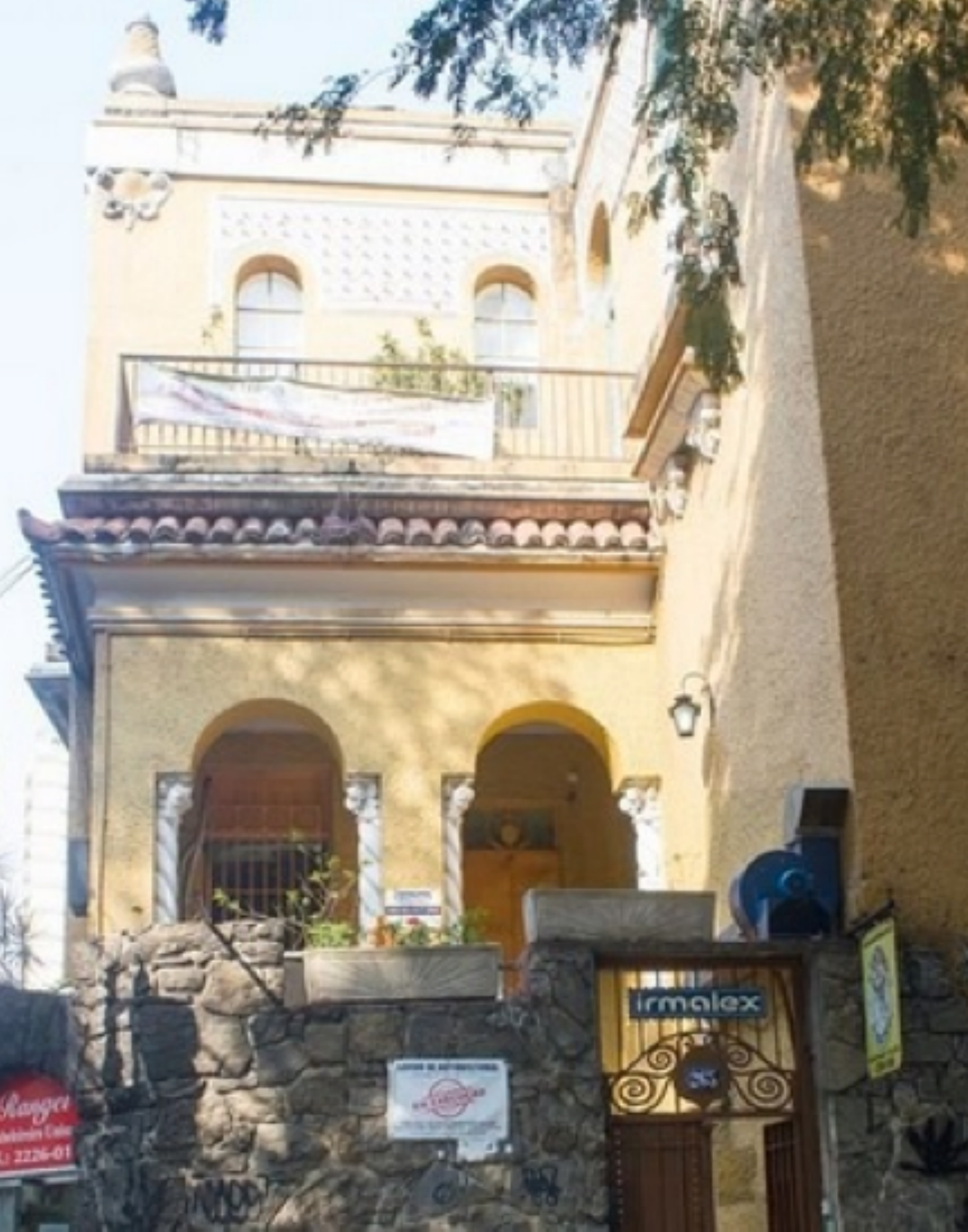}
\caption{\small Present view of the  facade  of the two store house,  located at the corner 
of the streets S\~ao Clemente and  Sorocaba, hosting IMPA in $1962 - 64.$  
 \label{fig:sorocaba}}
\end{center}
\end{figure}

The occasion was the Seminar   on the Qualitative Theory of  Differential Equations  (QTDE), directed by Prof.  Mauricio M. 
Peixoto (hereafter Prof.   
Peixoto, or simply Peixoto)  who had announced the following title  for his lecture: ``Open Problems on The Qualitative Theory of Differential Equations".
The seminar  activity had been interrupted  for a few  sessions.  
It was restarting after the return of Prof. Peixoto  who had travelled to attend  the International Congress of Mathematicians  at Stockholm, August  15 - 22.

\vglue .1in

 At that time Prof.  Peixoto  was the only Research Director resident  at IMPA.
 Prof. Leopoldo Nachbin (1922 - 1993), also a 
 research Director, was on leave of absence. 
 
 Peixoto  had very explicit
 views concerning Mathematics learning in graduate level. They where well known  to most participants:
``This science is assimilated through solving problems and thinking", he had clearly  stated, contesting insinuations
of those who favored  the preponderance of more intensive lecture courses and bookish  learning.
 This contrasted with the  naive and easy going
 vision I had partially acquired in my undergraduate predominantly reading contact with mathematics.

Excitement  and a tense expectation could be noticed in the audience.
For the most experienced participants of the Seminar,  the time to {\it face  true research problems}  had arrived. 
As a newcomer into advanced studies, and the youngest of all,  I could be  implicitly regarded outside  such group.

\section{ A flash 
  of
 mathematics 
at   IMPA  in  1962 - 63.}

\vglue .1in

Among the experienced  researchers, besides  Prof.  Peixoto   who stayed at IMPA the whole period from 
1962 to 
  1964,  were Elon  Lima (1930 - 2017),   Djairo  de  Figueiredo  (1934 - \;  ) and  Otto Endler (1929 - 1988)  who sojourned for  shorter periods.
Elon stayed the whole first semester and Djairo only a couple of months.
Otto stayed along 1962 and part of 1963.

Few lecturers  visited the Institute. Among them I mention below those which captured my interest.

Charles Pugh  (1940 - \; ). Subject: Closing Lemma, crucial in the work of Peixoto.

 Gilberto Loibel (1932 - 2013). Subject: Stratified Sets, introduced  by Ren\'e Thom  (1923 - 2004), founder, together with Hassler Whitney,  of the Theory of Singularities of Mappings. 

 Wilhelm Klingenberg (1924 - 2010). Subject: Closed Geodesics, located  in the intersection of Differential Geometry,  Differential Equations and Calculus of Variations.

At that time,  IMPA did not have a minimal program of courses to be  offered along the year.
This activity  depended on the research fellows present, in  a quickly changing regime.

However, there was a permanent basis providing  scientific stability to the Institute  along those years, around Prof. Peixoto research project.
Furthermore, due to an agreement with the University of Brazil, later denominated UFRJ (Federal University of Rio de Janeiro),  a doctoral program started at IMPA in 1963. 

\vskip 0.5 cm

{\bf  This pioneer research project, under the leadership of Peixoto, is the landmark of  the systematic interest  in Dynamical Systems (DS) in Brazil.   }

 {\bf This was the first explicit  
 effort to stimulate the initiation into  
 research in this area of mathematics in Brazil. 
 
 I had the unique chance and privilege to  be part of the first group of doctoral candidates, under
 the supervision of Prof. Peixoto.}

\vskip 0.5 cm
After the influential works of the American mathematician Stephen Smale (sixties and seventies) the denomination DS
nominally assimilated a significant part of the QTDE  from which the separating border is not well defined. See Smale's 1967 landmark article \cite{BAMS}.
DS is also the name of a famous  book of George Birkhoff,  printed in  1927.
\vskip .3 cm

Ordered by the preponderance they had in my initiation into research, below I list some of  the courses and seminars held at IMPA in 1962\footnote{This  abridged version concerns activities  on the QTDE. A wider  description can be found in \cite{l}.}:

1. Seminar on the Qualitative Theory of  Differential Equations;

2.  Seminar on Differentiable Manifolds following the inspiring Porto Alegre  Elon Lima  book  \cite{IVD};

3. Seminar based on the reading of the, now classic, book of J. Dieudonn\' e: Foundations of Modern Analysis \cite{dieudonne};

4.  Course on Algebraic Topology, taught by Peixoto, based on the books of Cairns and Hocking - Young;

5. Course on Multilinear Algebra and Exterior Differential Calculus, based on Bourbaki and  Flanders, taught by Elon Lima;

6. Seminar  on the reading of the book ``Lectures on Ordinary Differential Equations'' by W. Hurewicz \cite{hurewicz}.

I list bellow, with no ordering,  some  of the research fellows and {\it habitu\'ees}  of IMPA,  related to Prof. Peixoto Project, in 1962 - 63:

Ivan Kupka (1938 - \;   ),

 Maria L\' ucia Alvarenga  (1937 - \;   ),

Jacob Palis  (1940 - \;   ),

Lindolpho  C. Dias (1930 - \;   ),

Alcil\' ea Augusto, (1937 - \;   ),

 Aristides  C.  Barreto, (1935 - 2000).

Much  of the contents of  the activities  $1$ to  $5$  listed above   was  new  to me.
I devoted to them special attention,    with library search and extensive 
 complementary  readings.    In this endeavor the  friendly interaction with Ivan Kupka, by far the most knowledgeable in the group, was  auspicious.  

\vspace{-0.0cm}
\begin{figure}[h] 
\subfigure[
]
{ \includegraphics[scale=0.42]{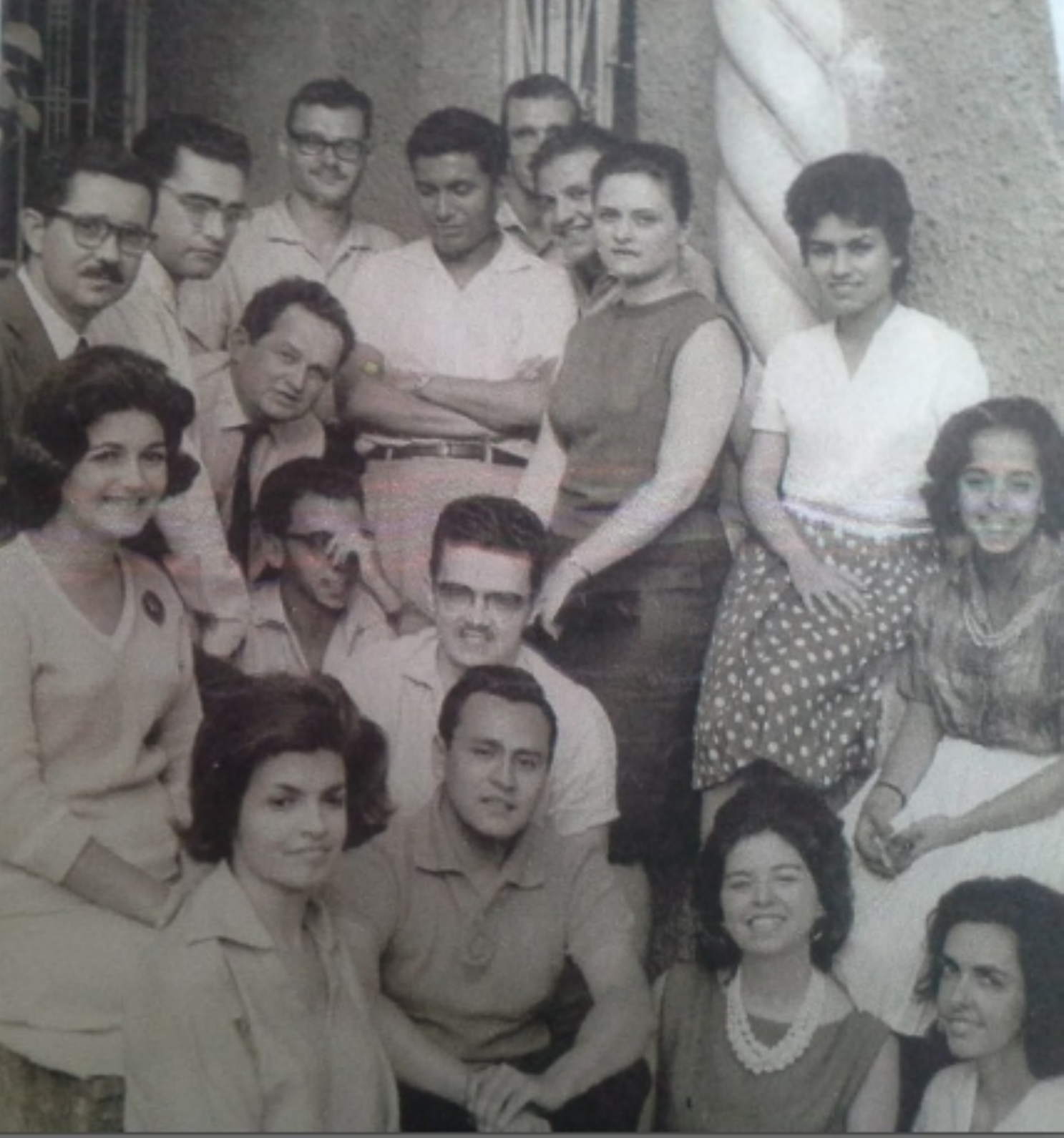}} 
\subfigure[
]
{ \includegraphics[scale=0.42]{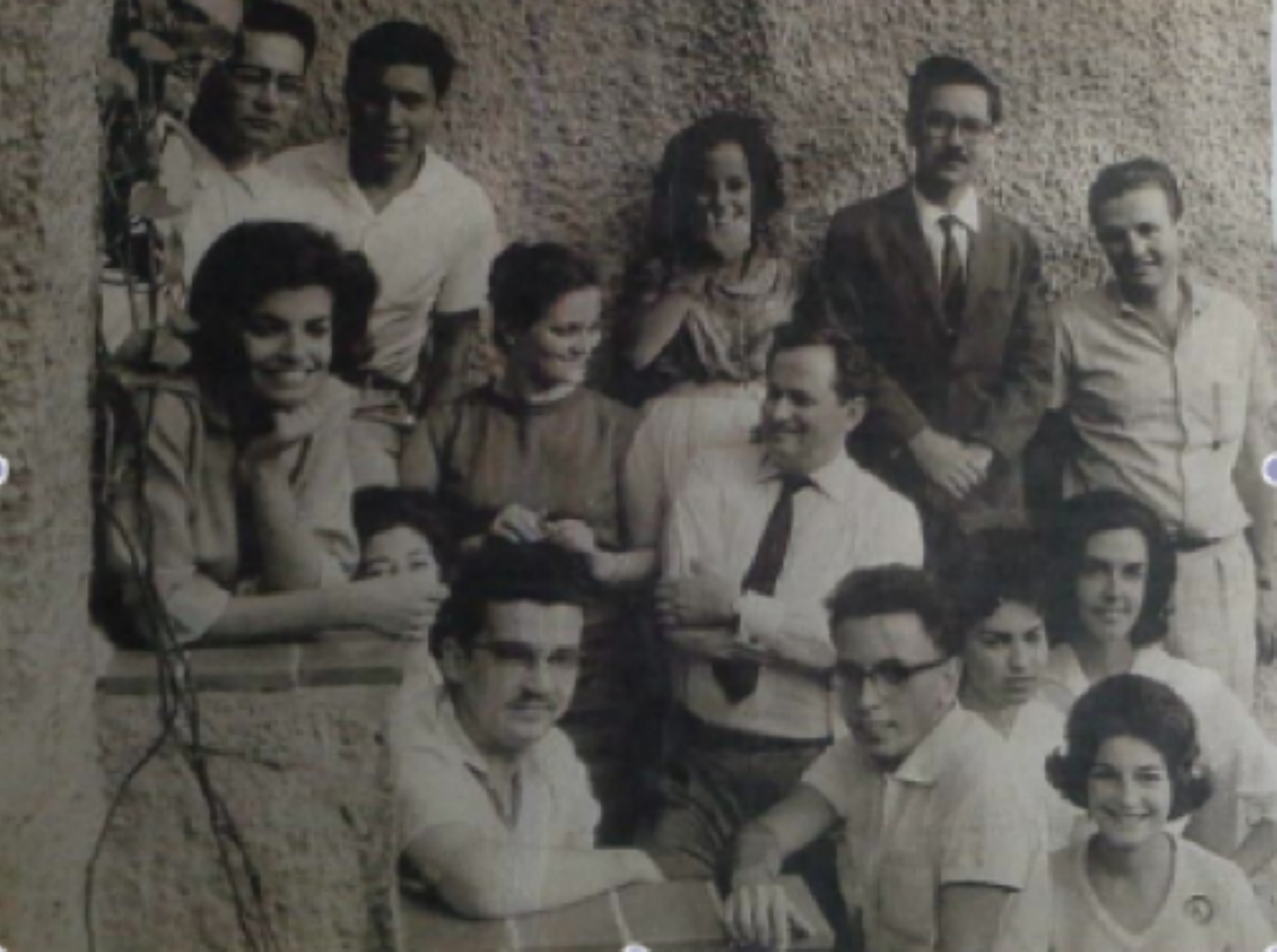}} 
\vspace{-0.0cm}
\caption{{\bf At the entrance of IMPA,  August 1962.} \newline
{\small  (a)  Top row, left to right: Lindolpho Dias,  M. Peixoto,   Ivan  Kupka, E. Isla, Aristides Barreto, Jacob Palis, Eliana Rocha, Alcil\'ea Augusto; descending the right margin  and continuing clockwise on the bottom margin:  C. Marquez, M.  L. Alvarenga, M. H. Cerqueira, J. Sotomayor, Lia Velloso and Adarcy P. Costa.  Three gentlemen sitting at the center, descending counter clockwise: L. Nachbin, H. Machado and J.  A. Barroso. \newline
(b) Top row from left to right:  M. Peixoto,  E. Isla, Celina Marquez,  Lindolpho Dias and Jacob  Palis. L. Nachbin at center, with necktie. Bottom right corner: Alcil\'ea   Augusto, M. L. Alvarenga,   Adarcy P. Costa. 
}}\label{Fig:02}
\end{figure}

\section {Peixoto's  Seminar on   QTDE, 1962.} \label{seminar}
Along  my  
sojourn
 at IMPA, 1962-64,  the  most  remarkable of all  the activities  
  in  which  I engaged   was the 
Seminar on the Qualitative Theory of Ordinary Differential Equations.
  
Following the epistolary reading directions sent to me by Peixoto  in 1961, recounted  in the evocative essay ``Mathematical Encounters"  \cite{encounters},  with the  inspiring   master presentations of Coddington  and  Levinson \cite{cod_lev}
and Hurewicz \cite{hurewicz},   I had already been initiated into the 
  first steps of the Qualitative Theory. 
 The phase portrait, limit cycles and singular points local structures: saddles, nodes, foci and centers, were fiound there. A bright synthesis of these elements in the  Poincar\' e - Bendixson Theorem complemented the introductory contact.
 
While the subjects presented and discussed  in the Seminar  did not require a   heavy  background knowledge, 
their in depth appreciation depended on
a level of maturity and on an inquisitive disposition, beyond the initiation outlined above,  to which I had not  yet been  exposed. 

\subsection {Non singular ODEs on the Klein Bottle and the Torus.}
The first lecture at the Seminar  TQEDO, was delivered by Peixoto at the beginning of April. 
The subject was  {\it The  Theorem of  Kneser, }
whose conclusion is that   ``Every vector field with no singular points  in the  Klein Bottle  has a periodic orbit''. See \cite{kneser}.

At the two next meetings  of the seminar,  Peixoto presented the basic theory of the {\it rotation number} following the last chapter of Coddington and Levinson.

The second presentation included also a geometric construction of  the {\it example of Denoy}
of  a $C^1$ non singular vector field 
for which all its orbits cluster  in a closed  invariant set, which transversally is a Cantor set, with no proper  subset sharing  these properties. 

Such a set is denominated  ``minimal non-trivial''. 
The ``trivial minimal''  sets 
are the  singular points --equlilibria-- and periodic orbits.  

\subsection {Invariant Manifolds.}
 Elon  Lima  continued the Seminar. He  presented  part of  Chap. 13 of Coddington and Levinson \cite{cod_lev}  
which contains  the Theory of Invariant Manifolds.
It deals with the n-dimensional generalizations of ``saddles''  for singular points --equilibria--  and periodic orbits. Actually,  these dynamical objects are called ``hyperbolic'', a name coined by Smale. This is  one of the most technical matters of the book, which follows the approach of the German mathematician Oskar Perron (1880 --1975).

In 1970,  C. Pugh, whose name will appear later in this essay, 
and M. Hirsh made a substantial extension of the Invariant Manifold Theory elaborating ideas in the works of the French mathematician Jacques Hadamard  (1865 - 1963). See \cite{hirsh_pugh}.
Years later, I had a better assimilation of \cite{hirsh_pugh} and, in 1979, included it in  \cite{licoes}.  In 1973, I had used it  to give a conceptual proof  of the smoothness of the flow of a vector field \cite{soto_bsbm_1973}.

Meanwhile,  the approach of Perron, also in  \cite{cod_lev},  was elaborated by C. Irwen  using the Implicit Function Theorem in Banach Spaces.  In \cite{cod_lev}   is used  the method of successive  approximations. 
A version of this idea can be found also in  Melo and Palis \cite{pal_mel}. 

\subsection {A detailed presentation of  three   papers of Peixoto on Structural Stability.} \label{AA_62}
Along part of May and June,  Alcil\' ea  Augusto delivered a series of presentations, 
very detailed and carefully  prepared, covering Peixoto's papers: \cite{peixoto_a}, \cite{peixoto_b}  and \cite{peixoto_t}. 
The last paper, however, involved  difficulties, particularly on non-orientable surfaces.  This 
matter evolved into the so called  {\it  The Closing Lemma Problem,}  
of present research interest,  extrapolating  the domain  of   Classical Analysis.  

Rather than outlining  the individual contents of these  papers, I  include  below a personal appreciation on the subject, with  non - exhaustive  references.

\section{ A Glimpse into Structural Stability.} 
The concept of Structural Stability  was established during the collaboration
of the 
Russian 
mathematicians A. Andronov ( 1901 - 1952) and L. Pontryagin   (1908 -1988)  that
started in 1932 \cite{pontryagin}. It first appeared in their research note published in 1937. 
Andronov (who was also a physicist) 
founded  
the  very  important  
Gorkii 
School of  Dynamical Systems. 
He left a
remarkable mathematical heritage, 
highly respected both in Russia and in the West, \cite{gorelik}.   
By  1932 Pontryagin was an already famous Topologist who had started  to  teach differential
 equations and had voiced his interest in studying  applied problems.

Structural Stability is a consequence of the encounter of two mathematical  cultures, See discussion in Sec.\ref{con_com}.  

For a dynamic model --that is, a differential equation or system 
$x^{\prime}=f(x)$-- to faithfully represent a phenomenon of the physical
world, it must have a certain degree of stability. Small perturbations,
unavoidable in the recording of data and experimentation, should not affect
its essential features. Mathematically this is expressed by the requirement
that the \emph{phase portrait} of the model, 
which is the geometric synthesis of the system, 
must be topologically unchanged by small perturbations. 
In other words, the phase
portraits of $f$ and $f+\Delta f$ must agree up to a homeomorphism of the
form $I+\Delta I$, where $I$ is the identity   
transformation of the phase space of the system and $||\Delta I||$ is small. 
A homeomorphism of the form $I+\Delta I$ is called an
$\epsilon$-{\emph homeomorphism}
if $||\Delta I||<\epsilon$; that is, it moves points
at most $\epsilon$ units from their original positions.

Andronov and Pontryagin 
 stated a characterization of structurally stable
systems on a disk in the plane. 
This work was supported by the analysis of numerous
concrete models of mechanical systems and electrical circuits, performed by Andronov and his associates \cite{a_L},  \cite{a}. The
concept of structural stability, initially called \emph{robustness}, 
represents a remarkable evolution of the 
continuation method of Poincar\'{e}. 

When the 
American mathematician S. Lefschetz translated the
writings of Andronov and his collaborators from Russian to English \cite{a_L}, 
he changed the name of the concept to the more
descriptive one it has today  \cite{a_L}.  He also stimulated 
H. B. de Baggis
to work on a proof of the main result as stated by Andronov and
Pontryagin.

Peixoto improved the results of the Russian pioneers in several directions. 

For example, he introduced the space 
$\mathcal{X}^r$ 
of all vector fields of class $C^r$ , 
and 
established the 
openness and genericity of  structurally stable vector fields on the 
plane and on orientable surfaces. 
He also removed the $\epsilon$-homeomorphism  
requirement from the original definition, proving that it is equivalent to the  existence of any homeomorphism. This was a substantial improvement of 
the Andronov-Pontryagin planar theory.

The transition from the plane to surfaces, as  
in 
Peixoto's work,  
takes us from classical ODEs to  the modern theory
of Dynamical Systems,   from Andronov and Pontryagin to D. V. Anosov and S. Smale. 

It has also raised 
delicate problems --for instance, the 
\emph{closing lemma}--  that have challenged mathematicians 
for decades \cite{gutierrez_closing}.

In \cite{sm} S. Smale regards Pei\-xo\-to's structural stability theorem as the 
prototypical example and fundamental model to follow for 
global analysis.

\section{
 Open Problems in ODEs, September 1962.  \label{list_english}} 

\vspace{-0.0cm}
\begin{figure}[!ht]
\begin{center}
\includegraphics[scale=0.99]{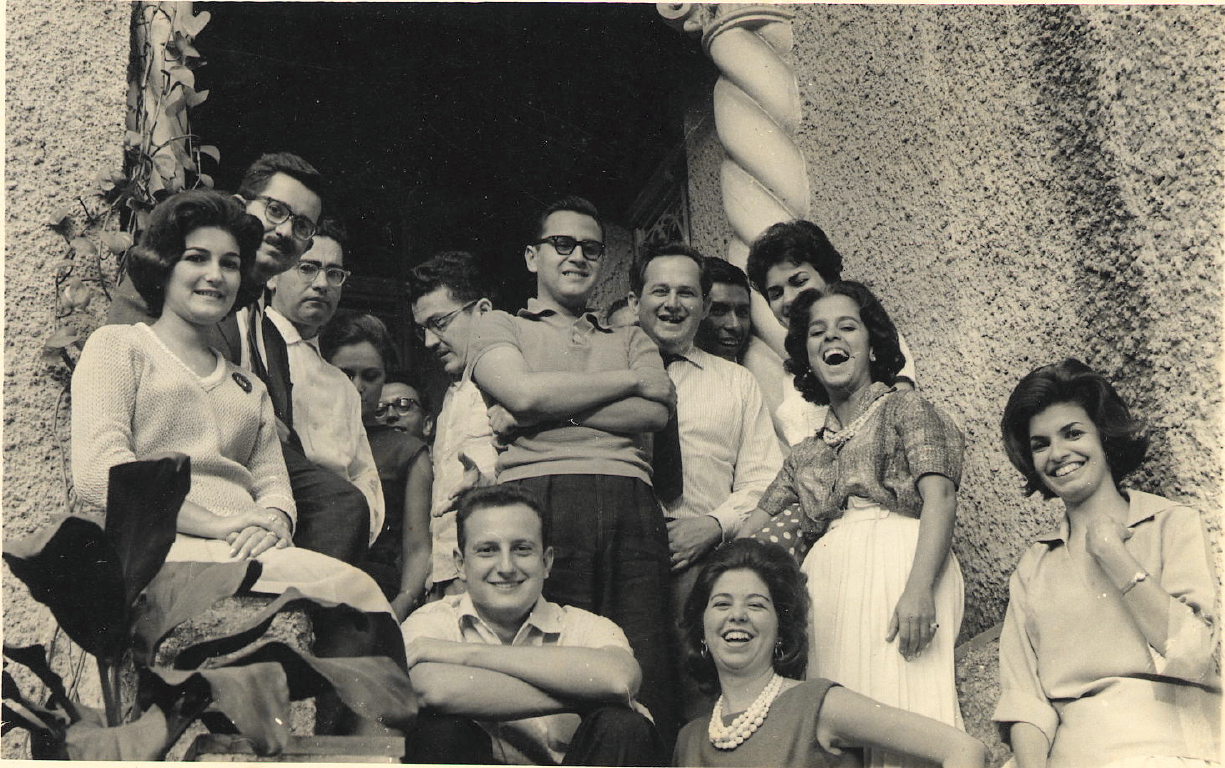}
\vspace{-0.0cm}
\caption{{\bf At the entrance of IMPA,  August 1962.} \newline
 {\small Standing, left to right facing camera: Lindolpho C. Dias, Mauricio Peixoto, Jorge Sotomayor, Leopoldo Nachbin,  
Alcil\'ea Augusto and  C. Marquez. Bottom row, left to right: Adarcy. P. Costa, Roberto  R. Baldino, M. H. Cerqueira and  Lia Velloso (librarian).  
}}\label{Fig:01}
\end{center}
\end{figure}

After a  concise, though  very emphatic,  introduction  about the importance  of attacking 
research problems, Prof. Peixoto began to enumerate and discuss  five of them.   
  
\subsection{First order structurally stable systems.} \label{p1}
\emph{Consider the complement
$\mathcal{X}^r _1$ of the set
$\Sigma^r$ 
of \ $C^{r}$-structurally stable
vector fields, relative to the set 
$\mathcal{X}^r$ of all vector fields  on a compact  two-dimensional manifold. Let 
$\mathcal{X}^r _1$ 
be endowed with the induced $C^{r}$ topology.
Characterize the set
$\Sigma^r_1$ of those vector fields that are
structurally stable with respect to arbitrarily small perturbations inside 
$\mathcal{X}^r _1$.  }

\smallskip

This problem goes back to a 1938 research announcement of A. A. Andronov and 
E. A. Leontovich 
\cite{al_38}, \cite{al}. 
They 
formulated 
a characterization of 
$\Sigma^r _1$ 
for a compact region  in the plane. 
This step points toward a
systematic study of the bifurcations (qualitative changes)
that occur in families of vector fields as they cross 
$\mathcal{X}^r _1$.
In the research announcement  --contained in a dense four pages note--  they 
stated  that the most stable
bifurcations 
occur in $\Sigma^r_{1}$,   \cite{s}, \cite{cbm_81}.

\subsection{The problem of the arc.} \label{p2}
\emph{ Prove or disprove that a continuous curve (an arc)
in the space $\mathcal{X}^r$ of vector fields of class $C^{r}$ on the sphere can
be arbitrarily well-approximated by a continuous curve that meets only 
finitely many \emph{bifurcation points}; 
that is, points  outside  the set of 
structurally stable vector fields, at which qualitative changes occur.}

Later research
established that 
 $\mathcal{X}^r_1$
enjoys great transversal complexity,
 which grows quickly with   the dimension of the phase domain. 
 
This  knowledge 
became 
apparent after the work of 
S. Smale \cite{BAMS} 
and also Newhouse  \cite{n}  and Palis -Takens  \cite{p_t}, among others.  

The  understanding of the phenomenon of persistent accumulation of bifurcations implies that  the problem of the arc as stated above has  a negative answer, \cite{IMUNI}. 
However, after removing the requirement of the approximation, 
Peixoto  
and
S. Newhouse proved that every pair of structurally stable vector fields
is connected by an arc that meets only finitely many bifurcation points. See \cite{p_n}.

\subsection{The classification problem.} \label{p3}
\emph{Use combinatorial invariants
to classify the connected components of the open set of structurally stable 
vector fields. 
}
\smallskip

The essential
difficulty of this problem is to determine when two structurally stable 
vector fields agree up to a homeomorphism that preserves their orbits and is isotopic to the
identity.

Some years later, Peixoto  himself  worked on  this  problem, \cite{peixoto_c}.

\subsection{The existence of nontrivial minimal sets.} \label{p4}
\emph{Do invariant perfect sets (that is, sets that are nonempty, compact, and 
transversally totally discontinuous) exist for differential equations of 
class $C^{2}$ on orientable two-dimensional manifolds?}

\smallskip
This problem goes back to H. Poincar\'{e} and A. Denjoy and was 
known to experts. 

It was solved in the negative direction by A. J. Schwartz \cite{sch_d}.  Peixoto presented this result
from a preprint 
that he received in November 1962. 

\subsection{Structurally stable second order differential equations.}\label{p5}
\emph{For equations of the form $x^{\prime\prime}=f(x,x^{\prime})$ 
(more precisely, for systems of the form $x^{\prime}=y,\,\,y^{\prime}=f(x,y)$), 
characterize structural stability, and prove the genericity of 
structural stability,  in the spirit  of  Peixoto's results for vector fields
on two-dimensional manifolds.}
 
 \vspace{0.3cm}
Problems \ref{p2} to \ref{p4}   were assigned, 
in one-to-one correspondence, to the
senior participants of  the seminar. 
The first and last  problems   were held in reserve for a few months.

In Sec. \ref{thesis},  I will recount how
  I was conduced to obtain 
Peixoto's support  to attack  problem \ref{p1}
on his list.

\vskip 0.3 cm

{\bf When,  years later,  I read  the proposal  of   the   famous  $1900$   
Hilbert Problems,  the
words in the epigraph
made me 
evoke     the  above mentioned   introit   of  Peixoto's  list of problems  in September 1962,  which  then,  keeping in mind the enormous  difference in proportions,     struck me
   as a    distant  echo of Hilbert's words  in Paris 1900, 
 that   reverberated along 
 decades 
 before reaching  the  tropics. }

\section {Peixoto's seminar, last 1962  sessions.}  \label{seminar_nov_62}

The last lecture in  the Seminar was delivered by Peixoto. 
 It was based on a preprint of   A. J. Schwartz \cite{sch_d},  which   solved  in the negative problem \ref{p4}.

Before this, toward the end of October,
Ivan Kupka  delivered  a series of very technical lectures about  
Invariant  Manifolds along arbitrary orbits, not necessarily periodic. 

He mentioned the work of Oskar Perron as the main reference for the hyperboliciy hypotheses adopted.

To mitigate the  concern of most participants in face  with abundant analytical technicalities involved in  the
rough presentations,  
Peixoto  started  a parallel series of  
very informal  tutorial discussions  on  Kupka's
 lectures. 

Everybody freely  expressed their disparate attempts  to explain  geometrically the ideas as well as the long  chains of inequalities involved.

Maybe  this convivial contact and my  
participation   in the discussions,  freely   formulating   hunches, 
established   a  more  direct channel of communication  between me and   Peixoto,  which so far had been a very formal  one.

\section{A   good research problem.\label{good_problem}}
After the presentations of 
Peixoto's work, timely  commented    by  himself,  complemented with  considerable  struggle with the bibliography,  
 it was possible  to have  a panoramic view of a fascinating  piece of  
  knowledge. 
  It was  a sample  of the evolution of mathematical  ideas  with an
  intriguing 
  historical background, mathematically deep but 
 essentially 
  accessible,   whose  consolidation had seen the light 
 in the last  four years,  concomitantly with my 
learning at  university level, of  the principles of   Mathematical Analysis, Geometry and Differential Equations.  

Complementing the seminar lectures  I had made substantial readings related to the QTDE and the Calculus and Geometry on Banach Spaces.
I devoted  considerable attention to Sard Theorem, which I first read in the  presentation,  \cite{EDJ},  of Prof. Edson  J\' udice of  the University of Minas Gerais (U.F.M.G), donated to me by the author by a recommendation of Aristides Barreto,  with whom I became close friend. 
I also studied the book by S. Lang  on Differentiable Manifolds, modeled on Banach spaces \cite{lang}.

Along this endeavor  I had rewarding discussions   on  mathematical subjects with Ivan Kupka.

With great  profit    I  read the 
enlightening 
lecture  notes   Introdu\c c\~ ao  \` a  Topologia Diferencial,  \cite{ITD} by  Elon Lima, which anticipated in several years Milnor's ``Topology from a differentiable viewpoint" \cite{MILNOR}. 
After this, I studied substantial parts of  L. S.  Pontrjagin, \cite{pontryagin_55},  ``Smooth Manifolds and its Applications to Homotopy Theory'',  which I consider  the original source for the application of Differential Analysis to the  study problems in Topology.

Some undefined intuitive and esthetic considerations,  and a certain  overestimation 
of  my   readings  about  Sard's Theorem,  
led me to the 
hunch 
that,
if not all,  part of Peixoto's  genericity Theorem of Structurally Stable  Systems  could be obtained from an appropriate 
infinite dimensional version of this theorem.
Being  undefined the involved domain and range spaces.

In discussions with Kupka,  I had learned of an extension in this direction due to Smale.

I approached Peixoto and  shared  with him my naive expectations. He made no comments. 
However,  toward  the end of  november he handed me a copy of the Andronov - Leontovich four page note \cite{al_38}. \\
He said: ``Do not loose this. It is important.  It is  a good problem".
 
 Despite the technical difficulties,  enhanced by, at  that time,     lack of bibliography on  bifurcations,  the  sensation of possessing  a research problem, produced in me the mixed feelings  of a naive fulfillment  and  of overwhelming responsibility.  
 
 I left Rio de Janeiro for the extended summer break of 1963. 
In my mind  I  carried   a   new sense of  mathematical awareness. 
In my bag,  
packed in a plastic cover, travelled  the note of Andronov - Leontovich. 

\begin{rem}
Concerning the 
 desideratum  of  providing a proof of Structural Stability Theorem genericity  theorem with Sard's Theorem, I mention that  
in 
the decade of 1970, when I began to lecture on the subject   for wide audiences,  I felt the 
need, and found,  a direct,  self-contained,   transparent
proof 
that worked for the plane and for polynomial vector fields. 
The course in the 1981, $13^{th}$
 Brazilian  Mathematics Colloquium,  \cite{cbm_81},  was an opportunity to communicate  the new  proof that used  only an elementary form  of the one-dimensional Sard's Theorem.  An abridged  version  was also published in \cite{cras_1981}.
\end{rem}

\section{ Back to IMPA in 1963.}  \label{back_63}

In December 1962  I took  the final examinations of  a   few courses I had pending  to fulfill  the number  of academic  credits  required  to get the  Mathematics Bachelor degree from the National University of San Marcos, Lima,  Peru.  

Along January and February I worked as teaching assistant in  a Summer   Mathematics School   for the training of  high school teachers.
In my spare time, and full time on March,  I scribbled piles of   pages of calculations and drawings  attempting  to decipher the statements in   Andronov- Leontovich  note. 

Arriving   to IMPA  at the end of March, I had some time to discuss with Ivan Kupka the outcome  of my summer struggle with  the note of Andronov - Leontovich. 
Prof. Peixoto  had scheduled   me  to present a report on the subject in May.  

The possibility of  adapting   the  methods of 
his works to get some form of density of 
${\Sigma}_1^r$ in ${\mathcal X}_1^r$, as proposed in problem \ref{p1},   was raised along 
 the preparation and in discussions after the seminar sessions.  This point is  not  
 mentioned in \cite{al_38} and  \cite{al}.

 At the expense of considerable work, it seemed 
possible to adapt the methods in the works of  Peixoto  to the   
 formulation extracted from the note of    Andronov - Leontovich, extending it   form the plane to surfaces.
 
 However,  while  attempting  to complement my personal studies  by    
   reading whatever work  containing 
 some material   on  bifurcations  that fell in my hands, among which ware  \cite{minorski} and  \cite{sansone-conti},  I was being led 
 to the suspicion   
 that   problems \ref{p1} and  \ref{p2}  were intimately interconnected. 

 This intuition
 received a mathematical formulation 
 during  July,  catalized   by events that took place  at  the  IV  Brazilian Mathematical Colloquium.

\section{The IV  Brazilian Mathematical Colloquium,  1963. \label{cbm}}    
   At the Fourth Brazilian Mathematics Colloquium,  in   July 1963, \cite{cbm_63}, \cite{l}, everyone  
 working  under Peixoto's supervision 
 made  reports 
 in a session of short communications. 
  Only  Ivan Kupka 
 delivered a plenary lecture.

Below, in free translation,  the titles  of  the communications:

Maria L\' ucia Alvarenga: ``Planar Structural Stability".   
 
Alcil\' ea Augusto:  ``Parametric Structural Stability",

 Aristides   Barreto:  ``Structural Stability of  Equations of the form  $ x^{\prime\prime}= f(x, x ')$".
 
  ``Higher order Structural Stability",  my communication,  closed the list.

 The part  of the work  I did  at IMPA  that I consider innovative,  concerning  the structure of smooth Banach sub-manifold of the 
class of Andronov - Leontovich vector fields, 
 had  not yet been conceived. 
However its  main ideas emerged during  the Colloquium.
I  explain this point  now.     

In the first semester of 1963,  Kupka had spent  a couple of months at the University of Columbia, N.Y., where Smale  worked.
He 
brought the, now classic,  paper of Palais and Smale   ``A  Morse Theory for Infinite dimensional manifolds". Encouraged by Peixoto, he presented a Plenary Lecture about it.  The title in  \cite{cbm_63}  is ``Counter example of Morse - Sard 
Theorem for  the case of Infinite Dimensional Manifolds", 
though  he spent most of his time explaining the Palais-Smale 
Theory.  
 
 After an enlightening presentation, he concluded:

``The introduction of  a new Theory  must be accompanied by a solid justification. 
The one I presented today has applications
to Calculus of Variations, Control Theory and Differential Geometry, among other subjects."  


In her short presentation Alcil\' ea Augusto  
exposed the generic finitude for the encounter 
of an arc of vector fields with  the class of  those whose equilibria  have  vanishing Jacobian.    
She used an  elementary version of Thom Transversality Theorem 
analogous to the one used to prove the invariance of the Euler -  Poincar\'e - Hopf 
 characteristic of  a smooth manifold, expressed as the sum of the indexes of the singularities of a vector field with non-vanishing Jacobians. 
At that time  I had seen this procedure  
clearly explained  in Lima's lecture notes 
\cite{ ITD}.

 Added to the examples of 
  one parameter bifurcations I had scribbled from  books on non-linear mechanics and oscillations, such as \cite{minorski}  \cite{a}, and the re-reading of \cite{al_38},  the ensemble of   the lectures  of A. Augusto and I.  Kupka,          impacted my vision of  problem \ref{p1}. 
 The  following intuition,  or   desideratum, 
struck
me:  The First and Second Problems were intimate  parts of the same problem. 
 The   suitable synthesis of this association 
  should be presented in terms of infinite dimensional sub-manifolds and transversality  to  sub-manifolds  in the Banach space of all vector fields tangent to a surface.  

In fact,  the  part of the  class ${\mathcal X}_1^r$  to be crossed  by an arc, or curve of vector fields, needed to  have a smooth structure,   
as that of a hypersurface, to express the crossing, or bifurcation, as a transversal intersection, thus unifying in one concept --the codimension--   all the diverse dynamical phenomena, including the global ones, such as the non-hyperbolic periodic orbits, homoclinic and heteroclinic orbits,  and not only the punctual ones, such as the singularities.    
  
The 
suitable approach to express this structure had to include  the mathematical  objects such as those  appearing in Kupka's lecture:
infinite dimensional manifolds with tangent spaces that could be used to express infinitesimally the transversal crossing with  ${\mathcal X}_1^r$.

The  expression  in  paper form
of this intuition had to wait some years  to see the light  \cite{s}.   

I will mention only two other lectures: \newline
- M.  Peixoto,  An elementary proof  of the Euler-Poincar\' e formula on Surfaces (``Uma prova elementar da formula de Euler-Poincar\' e em Superf\'\i cies''). 
\newline  
- Charles   Pugh,  ``The Closing Lemma''.

This Lemma is in fact an open problem whose statement for class  $C^s$ and dimension $n$ is as follows:  \newline 
{\it ``Every $C^r$ vector field on a compact n-dimensional manifold $M$  having a non trivial recurrent through $p \in M$ orbit  can be arbitrarily approximated  in the  $C^s , s \leq  r, $ topology by one which has through $p$  a periodic orbit."}

In his lecture for the 1963 Colloquium Pugh  presented a particular case for $n=2$ and $s=1$.   
Later he extended his analysis  for arbitrary $n$. 
Carlos Gutierrez made remarkable contributions to the case of $n=2, \, s>1$, where,  it  is still open, in most non- orientable surfaces \cite{gutierrez_closing}.  \\
The Closing Lemma Problem is a question that stems form Peixoto's works.

   {\bf Proportionally,  regarding  by subjects,  the  presence  of  Structural Stability and related topics --the school founded by  Peixoto-- in the  ensemble  of presentations in the  IV  Brazilian Mathematical Colloquium was overwhelming.}

\section{A thesis project based on the   good problem in  \ref{good_problem}.}  \label{thesis}

By  the end of November 1963,  the fellow 
that,   on  March 1962 had arrived to IMPA  with a bookish mathematical knowledge \cite{encounters}, seemed  somewhat 
 distant. 
 Prompted  
by the unfolding of Peixoto's Seminar, specially by the ODE Open Problems Session,  
outlined in Sec. \ref{list_english},  he  had  experienced  an upgrade 
on the amplitude of his mathematical imagery and  on the  profundity  and extension  of   his  knowledge. 
The readings performed and the mathematical events  at the IV Bras. Math. Colloq., 
outlined in Sect. \ref{cbm}  had a radical influence on his view
of the  problems of Peixoto presented  in Sec. \ref{list_english}.

\vspace{0.3cm}
It was clear then that  Problems $1$  and $2$  of  Peixoto's List were linked  by the differentiable  structure  of the extended class  of Andronov -  Leontovich that I denoted  ${\Sigma}_1^r$. 
The detailed analysis of this structure, however,  depended on making it explicit in several instances.

Peixoto agreed with 
 my  doctoral dissertation project  consisting  on the extension to surfaces of the class ${\Sigma}_1^r$, its  smooth structure  and  its density inside  ${\mathcal X}_1^r$.

An explicit counter example for the generic finitude of planar bifurcations was easy to find,  reconsidering in terms of transversality  and further elaborating the results concisely expressed by  the Russian pioneers.  
 See Sotomayor \cite{IMUNI}.

\section { The  1964 mathematical  works.} \label{works_64}

In  the  written composition  of  my  doctoral thesis  (DT)
was   deposited   most of  
what  I   had learned along  1962 - 63: \\
The Calculus in Banach spaces and manifolds,\\
The invariant manifolds {\it a la} Coddington-Levinson,\\
The Structural Stability papers of Peixoto,\\
The personal digest I had scribbled on the   understanding
of  Andronov - Leontovich (AL)
announcement note \cite{al}.

Keeping in mind the analogy with the previous works of Andronov and Pontrjagin, as improved by Peixoto, the note of AL, \cite{al_38},  can be outlined as consisting of:

1.- An axiomatic definition  for  the class   ${\Sigma }_1^r$  as the part of ${\mathcal X}_1^r$ =  ${\mathcal X}^r$${\backslash}$${\Sigma}^r$   that violate minimally the conditions of Andronov-
Pontrjagin and  Peixoto that  define  ${\Sigma}^r$.

 2.- A definition of the class  ${\mathcal S}_1^r$ of the systems in ${\mathcal X}_1^r$  that are structurally stable under   small perturbations inside ${\mathcal X}_1^r$.
 
 3.- The  statement  identifying  ${\Sigma }_1^r$ with  ${\mathcal S}_1^r$.

I transliterated  the terminology for the systems in \cite{al_38}  as  being  ``first  order structurally stable".  However Russian translators use the name : ``first order structurally unstable''.   In fact,  they are the  ``most stable among the unstable ones".   
See Andronov-Leontovich et al \cite{al}, where the proofs of the planar theorems of (AP) and (AL) were  published in 1971.

The programatic analogy between  Andronov - Leontovich  and  Andronov - Pontryagin, in planar domains,  is clear.  It strikes as natural to extend it to Peixoto's  surface  domains. However, DT takes this analogy further and prepares the way for a geometric synthesis  of the Generic Bifurcations.

In fact, it establishes the openness   and density  of ${\mathcal S}_1 ^r$  relative to the sub-space  ${\mathcal X}_1^r$. 
It also endows it  with the structure of 
smooth co-dimension one  sub-manifold 
  of    the Banach 
space ${\mathcal X}^r$.
 
This last, analytic and geometric,  aspect of DT has no parallel in the Russian works on the subject. 
It  makes possible to regard geometrically the simple bifurcations as the points of transversal intersection of a curve of  systems with the sub-manifold     ${\Sigma}_1^r$.

In DT are calculated the tangent spaces to each piece of ${\Sigma}_1^r$. 
For  the cases  of the homoclinic and heteroclinic  connections of saddle points, the functional whose kernel defines the pertinent tangent space, is expressed in terms of an improper convergent integral, which corresponds to the Melnikov Integral when restricted to vector fields,  i.e. autonomous systems.  In 1964 no reference  was  known to me.

 Sotomayor \cite{S} contains a study of the characterization of First Order Structural Stability in terms of Regularity of  ${\mathcal X}_1^r$.
 
\vskip 0.3cm
 The work of Ivan Kupka achieved celebrity after Peixoto published  \cite{peixoto_k}, which unified the versions of 
Smale, for diffeomorphisms,  and that of Kupka, for flows,  and coined the name {\it Kupka-Smale  Systems} for those  systems whose singularities  and periodic orbits are all hyperbolic and all pairs of associated  stable 
and unstable manifolds meet transversally.  

This work of Peixoto provided me with the language and  methods   that I had missed in 1964, for the extension of the class of Andronov Leontovich to a strictly larger immersed manifold, containing properly  the imbedded one, whose structure  was established in DT. 
Concerning   this immersed manifold, the transversality to it gives the generic position of an arc with a dense, in ${\mathcal X}_1 ^r$,  smooth part.  See  Sotomayor \cite{s}.  

\vskip 0.3cm

The work of Aristides Barreto studied  the systems of the form   $x' = y,   y' = f(x,y)$,  with  $f$  periodic in  $x$.
There,  he characterized  those which are structurally stable. 
As far as I know this is the first work on Structural Stability on non compact manifolds, the cylinder in this case. 
In \cite{s2} I presented a compact   version of the solution  problem \ref{p5}  in Peixoto's list.  
 
\section{Concluding  Comments. \label{con_com} }  
\subsection{ Timeline focused  in this 
essay, with some extrapolation.} 

 Looking in retrospect one   may be  tempted to think that some of the subjects presented  in the seminar, Sec. \ref{seminar},  had,  already in  1962-63,  some {\it scent of a distant past.}  However, it  cannot be denied  that they also glimpsed into the future.  In fact, for the next three or more  decades,  they had current interest  for an active line for  research training to work in Dynamical Systems, touching  its kernel.  Fundamental work  on these subjects was done along   the forthcoming years,  reaching   relatively recent ones. 
  On this matter allow me to evoke      
 the following universal words:
\begin{verse}
{\it If we wish to foresee the future of mathematics, \\our proper course is to study the history and\\  present condition of the science.
For us mathematicians, \\is not this procedure to some extent professional? \\We are accustomed to extrapolation, which is a method \\of deducing the future from the past and the present; \\and since we are well aware of its limitations, we run no risk  \\of deluding ourselves as to the scope of the results it gives us.}  
\\H. Poincar\' e, in The Future of Mathematics, \\read by G. Darboux in Rome, ICM, 1908.
\end{verse}

\vspace{0.5cm}
\begin{figure}[!ht]
\begin{center}
\vspace{-0.5cm}
\includegraphics[scale=0.50]{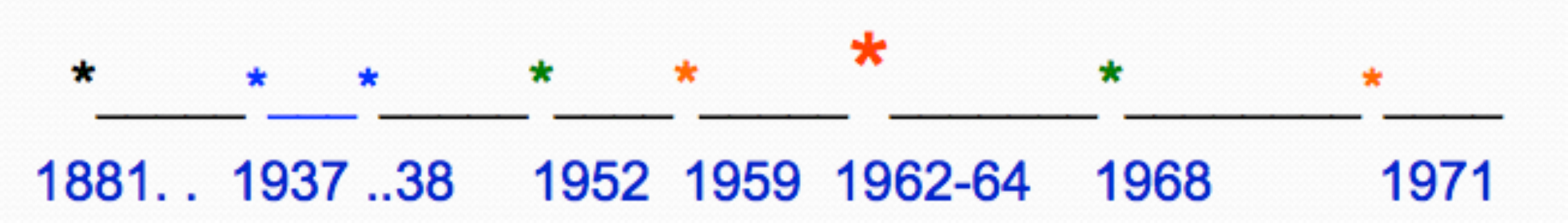}
\caption{{\bf Timeline with colored   landmarks,  weighted by the size of stars, 
with organizing center on the years  $1962 - 64$.} 
\newline
{\small Brazil (red): Peixoto' s works,  Seminar and Symposium;\newline
 France (black): Poincar\'e  QTDE;\newline
  Russia (blue): The Gorkii School Landmark; \newline
   USA (green): Lefschetz, 1949-52, and  Smale, Visit to IMPA, 1961, Seminar in Berkeley: 1966-67,  and his landmark papers Differentiable Dynamical System,  1967,  and What is Global Analysis?, 1969.
}} \label{Fig:02_1}
\end{center}
\end{figure}
\vspace{ -0.0cm}

\subsection{ Some Inquisitive   comments.} 
The mathematical concept of Structural Stability could hardly  have stemmed, in isolation, in  the offices of Mathematicians, pure or applied, or in the laboratories of Physicists or at the workshops of Engineers.
Something deeper and innovative happened in the collaboration of Andronov and Pontrjagin.

An effective collaboration involves the intellectual affinity 
of sprits. In this case,  involving  Andronov,  Physicist,   with exceptional mathematical knowledge, engaged in the research of the modeling of mechanisms, \cite{ gorelik},  
and Pontrjagin,  distinguished Mathematician, with remarkable contributions in Topology, interested in engaging himself in applied problems, \cite{pontryagin}.  

How do 
the transition from concrete examples and technological needs are processed  into seminal mathematical concepts and, afterwards,  to pertinent theorems?

  This is a
  central question of the Psychology  of the Creative Process, whose basis and analysis   
have been addressed by    \cite{hadamard} and  \cite{koestler}, among others.
\begin{verse}
``$\cdots$ {\it The creative act, by connecting  previously unrelated dimensions of experience,
 enables
 the authors to attain a higher level of mental evolution. 
 It is an act of liberation -- the defeat of habit by originality."} \\
{\small A. Koestler, (1964). The Act of Creation,  (p. 96). London: Hutchinson and Co.}
\end{verse}

However,  once formulated  in the domain of Mathematics, the concepts and theorems are amenable to generalizations, extensions and refinements, in style and essence.  Thus, they allow their  elaboration by Mathematicians, with their phantasies and the creative flight of their imagination.

 There are several stages in this transition in the realm of the evolution of mathematical ideas around Structural Stability, its extensions and generalizations. Maybe the first one, after the Russian pioneers, is that of Solomon Lefschetz, responsible for its  diffusion in the West  and  for  coining its expressive name, re-baptizing, the original  Robust Systems  given by the pioneers \cite{a_L}.  On this line of presentation, besides Peixoto, already cited, the names os Smale, Anosov, Arnold, Thom  and Mather, among others, should be mentioned, thus extrapolating  the realm of Differential Equations and Dynamical  Systems.

What, in an attempt of  expository simplification,   I referred  to above  as the outcome  of  the encounter of two distinct mathematical cultures:   
 knowledgable expertise and mathematical talent, \cite{gorelik}, \cite {pontryagin},  
 may,  perhaps,   be  better  
 explained in the delicate  threshold  between Mathematics and Art.

\begin{verse}
{\it Mathematics, rightly viewed, possesses not only truth,\\ but supreme beauty, a beauty cold and austere, \\ like that of sculpture, without appeal to any part of our weaker nature, \\without the gorgeous trappings of painting or music, yet sublimely pure, \\
and capable of a stern perfection such as only the greatest art can show.}\\
B. Russell, (1919)  "The Study of Mathematics", Mysticism and Logic \\ 
and  Other Essays. Longman. p. 60. 
\end{verse}

The mathematical and philosophical implications of Structural Stability,
can be appreciated in its  extensions to higher dimensional Dynamical Systems
 and to other domains of  the  Analysis on Manifolds,  such as the Singularities of Differentiable Mappings and the Theory of Catastrophes,  \cite{thom} \cite{sing_76},  and Multiparametric  Bifurcations, \cite{DRSe}  \cite{DRSz} \cite{cusp_loops}, as well as to  Classical Differential Geometry, such as the configurations of principal curvature lines and umbilic points on surfaces \cite{sg}, \cite{monge},  \cite{ga_so},  \cite{A_M}.

\vskip 0.5cm
\noindent {\bf Aknowledgement. }  Thanks are due 
 to  Lev Lerman for  sending me a copy of \cite{gorelik},  to Misha Dokuchaev for  helping  in its translation and   to Maria Lucia Alvarenga for 
kindly 
lending the photos in Figs.  \ref{Fig:02} and \ref{Fig:01}.

\vskip 1cm
 
\noindent Instituto de Matem\' atica e Estat\'\i stica,

\noindent
 Universidade de S\~ ao Paulo \newline
Rua do Mat\~ ao 110,  Cidade Universit\' aria,\newline
05508-090 \;     S\~ ao Paulo,  S.P., Brasil \newline
e-mail: sotp@ime.usp.br

\end{document}